\documentclass[a4paper]{amsart}

\usepackage{graphicx}
\usepackage{amsmath}
\usepackage{amssymb}
\input xy
\xyoption{all}

\renewcommand{\mod}{\operatorname{mod}\nolimits}

\newcommand{\End}{\operatorname{End}\nolimits}

\newcommand{\Ext}{\operatorname{Ext}\nolimits}
\newcommand{\Hom}{\operatorname{Hom}\nolimits}

\renewcommand{\dim}{\operatorname{dim}\nolimits}

\newcommand{\op}{\operatorname{op}\nolimits}
\newcommand{\Fac}{\operatorname{Fac}\nolimits}

\newcommand{\G}{\Gamma}
\renewcommand{\L}{\Lambda}
\newcommand{\C}{\operatorname{\mathcal C}\nolimits}

\newcommand{\ZZ}{\operatorname{\mathbb Z}\nolimits}

\newtheorem{thm}{Theorem}[section]

\newtheorem{algorithm}[thm]{Algorithm}
\newtheorem{cor}[thm]{Corollary}
\newtheorem{conj}[thm]{Conjecture}

\newtheorem{lem}[thm]{Lemma}
\newtheorem{prop}[thm]{Proposition}

\begin{document}
\title[Minimal infinite type]{Tame concealed algebras and cluster quivers of minimal infinite type}

\author[Buan]{Aslak Bakke Buan}
\address{Institutt for matematiske fag\\
Norges teknisk-naturvitenskapelige universitet\\
N-7491 Trondheim\\
Norway}
\email{aslakb@math.ntnu.no}

\author[Reiten]{Idun Reiten}
\address{Institutt for matematiske fag\\
Norges teknisk-naturvitenskapelige universitet\\
N-7491 Trondheim\\
Norway}
\email{idunr@math.ntnu.no}

\author[Seven]{Ahmet I. Seven{*}}
\thanks{(*) Supported by EU-LieGrits postdoctoral fellowship} 
\address{
Middle East Technical University\\
Ankara\\
06531\\ 
Turkey}
\email{aseven@metu.edu.tr}


\begin{abstract}
The well-known list of Happel-Vossieck of tame concealed algebras in terms of quivers with relations, 
and the list of A. Seven of minimal infinite cluster quivers are compared. There is a 1--1 correspondence
between the items in these lists, and we explain how an item in one list naturally corresponds to
an item in the other list. A central tool for understanding this correspondence is the
theory of cluster-tilted algebras.  
\end{abstract}

\maketitle

\section*{Introduction}

In the representation theory of finite dimensional algebras there 
is a famous list of algebras by Happel-Vossieck \cite{hv}
of the tame concealed algebras in terms of quivers with relations.
These constitute, together with the so-called generalized Kronecker algebras, 
the algebras $\L$ of minimal infinite type having a preprojective component
(see Section \ref{sec_back} for definitions). These algebras are useful for testing whether
finite dimensional algebras are of finite representation type. 

A. Seven recently produced another list \cite{se} in connection with his work 
on cluster algebras.
Cluster algebras were defined and first studied by Fomin and Zelevinsky \cite{fz1}.
We consider the special case (the case of ``no coefficients'' and skew-symmetric matrices)
where these algebras are determined
by a  finite quiver $Q$ with no loops and no oriented cycles of length two. 
We call such a quiver $Q$ a {\em cluster quiver}. 
Here we work over an 
algebraically closed field $k$.

A central concept
in the theory of cluster algebras is the mutation of quivers, and the quivers in the same mutation
class determine the same cluster algebras. The Dynkin quivers
$A_n, D_n, E_6, E_7, E_8$ and their mutation classes correspond to cluster
algebras of finite type \cite{fz2}. 
The list produced in \cite{se} gives the underlying
graphs of quivers with at least three vertices with the following properties. 
\begin{itemize}
\item[-]{The quivers in the list are not mutation equivalent to a Dynkin quiver}
\item[-]{Whenever a vertex is removed from a quiver on the list, the resulting quiver is mutation
equivalent to a Dynkin quiver}
\end{itemize}
The first condition says that the cluster algebras defined by these quivers are not of finite
type, while the second condition says that a cluster algebra defined by
any quiver obtained by removing one vertex is of finite type. 
The quivers from A. Seven's list in addition to the  
generalized Kronecker quivers (which have two vertices) are the only quivers satisfying these conditions.
Such quivers are here called {\em minimal infinite} cluster quivers.

The quivers of Happel-Vossieck contain also some dotted arrows, which correspond
to relations. If one replaces the dotted arrows in their list with
solid ones in the opposite direction, one obtains exactly the simply-laced 
quivers on A. Sevens list \cite{se}.

The aim of this paper is to explain why these lists are so closely related. For this we use the theory 
of cluster categories and cluster-tilted algebras from \cite{bmrrt, bmr1, bmr2}.
This theory was motivated by trying to model the essential ingredients in the definition 
of a cluster algebra in module theoretical/categorical terms. In particular we explain why and 
how we can construct one list from the other one. 
A key point is that by using \cite{bmr2} we see that the quivers of minimal infinite
cluster algebras coincide with the quivers of the cluster-tilted algebras 
of infinite type where all factor algebras obtained by factoring out an ideal generated by a vertex
are of finite type.
Then we compare the conditions of minimal infinite type for tilted and cluster-tilted 
algebras, and give a procedure for passing from the first class to the second one,
as was done in \cite{br} for tilted algebras of Dynkin type.

For going in the opposite direction we need to remove arrows from the quivers in the list of A. Seven.
Actually the list of A. Seven is a list of graphs, and the associated 
quivers are those where every full cycle is an oriented cycle, except in the case of $\widetilde{A_n}$,
where we must have a non-oriented cycle.

Here we use quadratic forms of quivers with solid and dotted arrows to describe which set of 
solid arrows should be made dotted, i.e. removed from the quiver of a minimal infinite cluster-tilted algebra
in order to recover the quiver for the corresponding tame concealed algebra. We also provide an 
algorithm for actually finding our desired set of arrows.

The results of the first four sections were announced at the Oberwolfach
meeting ``Representation Theory of Finite-dimensional Algebras'' in February 2005, 
where we also thank Thomas Br\"{u}stle for interesting conversations.

Assem, Br\"{u}stle and Schiffler have generalized Proposition
\ref{prop4.1} to any finite dimensional tilted algebra, see \cite{abs}.

\section{General background}\label{sec_back}
In this section we give some relevant background material from the
representation theory of finite dimensional algebras, from the theory of cluster algebras, 
and from the theory of cluster categories and cluster-tilted
algebras. 

\subsection{Finite dimensional algebras}
Let $\Lambda$ be a connected finite dimensional algebra over an algebraically
closed field $k$. 
Then the category $\mod \L$ of finite dimensional (left) $\L$-modules has almost split 
sequences, and there is an associated translation $\tau$.
For some finite dimensional algebras the AR-quiver can have special components
called  \emph{preprojective} components. These components are defined by the following
property:
each indecomposable module in the component lies in the $\tau$-orbit 
of an indecomposable projective $\Lambda$-module, and the component has no
oriented cycles \cite{hr1}. We later use that only a finite number
of indecomposable $\Lambda$-modules, up to isomorphism, have nonzero
maps to a given $\Lambda$-module in a preprojective component \cite{hr1}.

Not all algebras have preprojective components. But the hereditary
algebras $H=kQ$, where $Q$ is a finite quiver without oriented cycles, have
such components. Also the tilted algebras have preprojective components
\cite{st}. Recall that an algebra $\Lambda$ is tilted if $\Lambda =
\End_{H}(T)^{\op}$ where $T$ is a tilting module over a hereditary
algebra $H$, that is, $\Ext_{H}^1(T,T) =0$ and there is an exact
sequence $0 \to H \to T_0 \to T_1 \to 0$ with $T_0$ and $T_1$ direct summands
of finite direct sums of copies of $T$. The quiver of a tilted algebra is known to have no
oriented cycles. A central class of tilted
algebras are the \emph{tame concealed algebras}, that is, the
algebras of the form $\End_{H}(T)^{\op}$, where $H=kQ$ for an extended
Dynkin quiver $Q$, and $T$ is a \emph{preprojective tilting
module}, that is, $T$ lies in the (unique) preprojective component of
$H$.

The quiver with two vertices $a,b$ and 
with $r \geq 2$ arrows from $a$ to $b$ is often called a {\em generalized Kronecker quiver},
and we call a path algebra of such a quiver a {\em generalized Kronecker algebra}.

\subsection{Cluster algebras}
We here indicate the definition of a cluster algebra. 
The general definition from \cite{fz1} involves ``coefficients'' and 
skew-symmetrizable integer matrices. We restrict to the case 
with no coefficients and skew-symmetric matrices. 

Start with a \emph{seed}
$(\underline{x},Q)$, where $\underline{x}=\{x_1, \cdots, x_n\}$ is a
transcendence basis for the rational function field $\mathbb{Q}(x_1,
\cdots ,x_n)$ and $Q$ is a finite quiver with no loops or cycles of
length two. For each $i=1, \cdots , n$ 
one can define a new quiver $\mu_i (Q)$
and a new seed $\mu_i ((\underline{x},Q))=
(\underline{x'},Q')$, where $\underline{x'}=\{x_1, \cdots ,x_i',
\cdots ,x_n\}$ is a transcendence basis where the element $x_i$ in
$\underline{x}$ is replaced by an element $x_i'$ in  $\mathbb{Q}(x_1,
\cdots ,x_n)$, depending on $\underline{x}$ and $Q$. 
The new quiver $Q'$ (or seed $(\underline{x'},Q')$) is called the {\em mutation} of 
the quiver $Q$ (or seed $(\underline{x},Q)$) in direction $i$.
Continuing this way with mutation of quivers and seeds, the {\em cluster algebra}
$\mathcal{A}(Q)$ is defined to be the subalgebra of  $\mathbb{Q}(x_1,
\cdots ,x_n)$ generated by all elements appearing in the $n$-element
subsets $\underline{x},\underline{x'}, \cdots$. These subsets are called
{\em clusters}. 
The quivers obtained by a sequence of mutations of the quiver $Q$ are said to be
mutation equivalent to $Q$.

A major result in \cite{fz2} is that there is only a finite
number of seeds if and only if the quivers occurring in seeds
are mutation equivalent to a Dynkin quiver.

\subsection{Cluster categories and cluster-tilted algebras}\label{sec1.3}
Let $H=kQ$ again be a finite dimensional hereditary $k$-algebra, and
let $\mathcal{D}^b(H)$ denote the bounded derived category of finitely
generated $H$-modules. Denote by $\tau \colon \mathcal{D}^b(H) \to
\mathcal{D}^b(H)$ the equivalence such that for $C$ an indecomposable
object in $\mathcal{D}^b(H)$ we have an almost split triangle $\tau C
\to B \to C \to $. We then also have an equivalence $F=\tau^{-1}[1]$,
where $[1]$ denotes the shift functor in $\mathcal{D}^b(H)$. The
cluster category $\mathcal{C}_{H}$, introduced and investigated in
\cite{bmrrt}, is by definition the factor category
$\mathcal{D}^b(H)/F$, which is a triangulated category \cite{k}. 

A central concept is the notion of (cluster-)tilting object $T$ in
$\mathcal{C}_{H}$, where \\ $\Ext^1_{\mathcal{C}_{H}}(T,T)=0$ and $T$ is
maximal with this property, that is, if 
$\Ext^1_{\mathcal{C}_{H}}(T\oplus X, T\oplus X)=0$, then $X$ is a
direct summand of a finite direct sum of copies of $T$. It is shown that any
tilting $H$-module induces a tilting object in $\mathcal{C}_{H}$, and
by possibly changing $H$ up to derived equivalence, all tilting
objects in $\mathcal{C}_{H}$ are obtained this way.

A study of the closely related cluster-tilted algebras $\Gamma
=\End_{\mathcal{C}_{H}}(T)^{\op}$, where $T$ is a tilting object in
$\mathcal{C}_{H}$, was initiated in \cite{bmr1}. 
A useful result is that if $\G$ is cluster-tilted and $e$ is a vertex in the quiver,
then also $\G/\G e \G$ is cluster-tilted \cite{bmr2}. As a consequence we have the technique of 
shortening of paths from \cite{bmr3}, namely if $a \longrightarrow b \longrightarrow c$
is a path in a quiver with non-zero composition and no other arrows
leaving or entering $b$, then we can replace this path by $a \longrightarrow c$, and 
the new
quiver is still the quiver of a cluster-tilted algebra. 
Assume the quiver $Q$ has only single arrows.
We say that a (not necessarily oriented) cycle is {\em full} if the subquiver generated by the cycle
contains no further arrows.
For an arrow $\alpha \colon j \to i$ in the quiver of a cluster-tilted algebra we say that
a path from $i$ to $j$ is a {\em shortest path} if it does not go through any oriented cycle,  
and together with $\alpha$ gives a full oriented cycle.
A shortest path from $i$ to $j$ is necessarily involved in a 
relation from $i$ to $j$.
For finite representation type we have that any full cycle is oriented, and that the homomorphism space
between two vertices (more precisely between the corresponding projective modules) is
at most one-dimensional. 
For finite type there are at most two shortest paths from 
$i$ to $j$ corresponding to $\alpha$. One shortest path gives rise to a zero-relation and two 
shortest paths give rise to a commutativity relation \cite{bmr3}.


\section{Preliminaries}
In this section we discuss the two classification theorems which we
compare in this paper.

We start with discussing the work of Happel-Vossieck \cite{hv}. A basic connected finite
dimensional algebra $\Lambda$ is said to be of minimal infinite type if
it is of infinite type and for each vertex $e$ in the quiver of
$\Lambda$ we have that $\Lambda / {\Lambda e \Lambda}$ is of finite
type. There is no general classification theorem of algebras of finite
type in the sense that we have a list of all of them. There is however
the following useful criterion for when an algebra is of minimal
infinite type. 

\begin{thm}[\cite{hv}]\label{thm2.1}
An algebra is
of minimal infinite type and has a preprojective component
if and only if it is a tame
concealed algebra or it is a generalized Kronecker algebra.
\end{thm}

Actually a list of all basic connected tame concealed algebras 
in terms of quivers with relations
is provided in 
\cite{hv}. The above result has been even more useful because of this list.

Since we have seen that the tilted algebras always have a
preprojective component, a special case of Theorem \ref{thm2.1} can be
formulated as follows.

\begin{cor}\label{cor2.2}
The minimal infinite type tilted algebras 
are the tame concealed algebras and the generalized Kronecker algebras.
\end{cor}

As we have seen, cluster algebras (with ''no coefficients'') are
determined by finite quivers with no loops or cycles of length
two, i.e. cluster quivers. 
The quivers associated with cluster algebras of finite type are
the mutation (equivalence) classes of the Dynkin quivers. 
There is also here no known
list of the class of quivers obtained this way. But A. Seven \cite{se}
has given a list of graphs which can be made into quivers by choosing
a direction of arrows such that each full cycle becomes oriented. 
There is one exception, if the underlying graph is $\widetilde{A_n}$,
then the possible orientations are the ones which do not give an oriented
cycle. 
The associated list
of finite quivers contains exactly the
finite quivers (called minimal infinite cluster quivers)
with the property that if any vertex is removed, then
we get into the mutation class of a Dynkin diagram.
Note that it is a consequence of A. Seven's result that all 
minimal infinite cluster quivers are mutation equivalent 
to a quiver with no oriented cycles.

It was observed in \cite{se} that the lists are very closely
related, and they give in fact the same quivers if the
dotted edges in the Happel-Vossieck list (which we can view as dotted arrows) 
are replaced by solid
arrows in the opposite direction. 
In the remaining part of this paper we will give an explanation for why
this is the case, using the theory of cluster-tilted algebras. We will 
also give a procedure for how to construct one list from the other
one, and explain why it works.

\section{Interplay}
In this section we discuss the theory behind the relationship between
the two lists.
For this it is useful to note the following interpretation, using \cite{bmr1,bmr2}, 
of the quivers appearing in A. Seven's list. 

\begin{thm}\label{thm3.1}
Let $Q$ be a cluster quiver.
Then $Q$ is minimal infinite if and only if $Q$ is the quiver of a
basic cluster-tilted algebra $\Lambda$ of infinite type with the property that
$\Lambda/{\Lambda e \Lambda}$ is of finite type for each vertex $e$ in
the quiver.
\end{thm}

\begin{proof}
In \cite{bmr2} it was shown that for a finite connected quiver $Q$
with no oriented cycles, the quivers in the mutation class of $Q$
coincide with the quivers of the cluster-tilted algebras 
coming from the cluster category $\mathcal{C}_{H}$ for $H=kQ$. We also know that a
connected cluster-tilted algebra is of finite type if
and only if the corresponding algebra $H=kQ$ is of finite type \cite{bmr1}, 
that
is, if and only if $Q$ is a Dynkin quiver. 
Recall also from Section \ref{sec1.3} that for any cluster-tilted algebra
$\Gamma$ with a vertex $e$ in the quiver of $\Gamma$ we have that
$\Gamma /{\Gamma e \Gamma}$ is again a cluster-tilted
algebra.

Assume $Q$ is a minimal infinite cluster quiver. Then
it follows that
a cluster-tilted algebra $\G$ with quiver $Q$ is of infinite type
and whenever a vertex is removed from $Q$ the reduced quiver $Q'$
is the quiver of a cluster-tilted algebra of finite type.
Hence $Q'$ is in the 
mutation class of a disjoint union of Dynkin quivers. 
It is clear that also the converse holds.
This finishes the
proof. 
\end{proof}

While tilted and cluster-tilted algebras are quite different with
respect to homological properties, they are closely connected through
the fact that they are constructed from a common tilting
module. Theorem \ref{thm3.1} shows that the 
list of A. Seven gives a list of cluster-tilted algebras 
with properties similar to 
properties of the tilted algebras appearing in
the Happel-Vossieck list. This motivates the following.

\begin{thm}\label{thm3.2}
Let $T$ be a tilting module over a finite dimensional hereditary
$k$-algebra $H=kQ$, and let $\Lambda$ be the tilted algebra
$\End_{H}(T)^{\op}$ and $\Gamma$ the cluster-tilted algebra
$\End_{\mathcal{C}_{H}}(T)^{\op}$. Assume that $\Lambda$ is of
infinite type. Then $\Lambda$ is of minimal infinite type if and only
if $\Gamma$ is of minimal infinite type.
\end{thm}

\begin{proof}
This is obvious if $\L $ has at most two simple modules.
Assume first that the tilted algebra $\Lambda$ is of  minimal infinite
type, so that $\Lambda$ is tame concealed. Since $\Lambda$ is a factor
algebra of $\Gamma$, it follows that $\Gamma$ is also of infinite
type. Consider the diagram

$$\xymatrix{
\mod H \ar[rr]^{\Hom_{H}(T,\ )} \ar@{^{(}->}[d] & & \mod \Lambda
  \ar@{^{(}->}[d]^i \\
 \mathcal{C}_{H}  \ar[rr]^{\Hom_{\mathcal{C}_{H}}(T,\ )} & & 
  \mod \Gamma 
}$$ 
where $i$ is the natural inclusion functor obtained from
$\Lambda$ being a factor algebra of $\Gamma$. 
Since $\Lambda$ is tame concealed, we can assume that the tilting
module $T$ is preprojective. Then all but a finite number of
indecomposable $H$-modules are in the subcategory $\Fac T$ of $\mod
\Lambda$, whose objects are the factors of finite direct sums of
copies of $T$. Denote by $(\Fac T)_0$ the subcategory of $\Fac T$
where the $H$-modules $X$ in $(\Fac T)_0$ have the property that
$\Ext^1(T, \tau^{-1}X) \simeq D \Hom (X, \tau^2 T) = 0$.
Also all but a finite number of indecomposable $H$-modules are in
$(\Fac T)_0$, and the same holds when considering $(\Fac T)_0$ as a
subcategory of $\mathcal{C}_{H}$. Since $\Hom_{D^b(H)}(T, \tau^{-1}X[1])=0$
for $X$ in $(\Fac T)_0$, it follows that 
$\Hom_{\mathcal{\C}_{H}}(T, \ )\mid_{(\Fac T)_0}$ has its image in $\mod \Lambda$, 
and is the same 
subcategory as the image of $\Hom_{H}(T, { })\mid_{(\Fac T)_0}$. Since
only a finite number of indecomposable $\Gamma$-modules are not in
this image, it follows that only a finite number of indecomposable
$\Gamma$-modules are not $\Lambda$-modules.
In particular, for each vertex $e$ in the quiver of $\Gamma$, there is
only a finite number of indecomposable $\Gamma /{\Gamma e
\Gamma}$-modules which are not $\Lambda /{\Lambda e
\Lambda}$-modules. Since $\Lambda /{\Lambda e \Lambda}$ is of finite
type, it follows that $\Gamma /{\Gamma e \Gamma}$ is of finite type.

Assume now that the cluster-tilted algebra $\Gamma$ is of minimal
infinite type. The corresponding tilted algebra $\Lambda$ is of
infinite type by assumption, and $\Lambda /{\Lambda e\Lambda}$ is of
finite type since it is a factor algebra of the finite type algebra
$\Gamma /{\Gamma e \Gamma}$.
\end{proof}
\textbf{Remark:} The assumption that the tilted algebra $\Lambda$ is
of infinite type can not be dropped. For if a tilting module $T$ over
a tame hereditary algebra $H$ has both a nonzero preprojective and a
nonzero preinjective direct summand, then $\Lambda =
\End_{H}(T)^{\op}$ is of finite type \cite{hr2}, while the
cluster-tilted algebra $\Gamma = \End_{\mathcal{C}_{H}}(T)^{\op}$ is
of infinite type since $H$ is of infinite type.

\section{From tilted to cluster-tilted algebras}
In this section we use the previous results to show how and why the
list of A. Seven can be obtained from the Happel-Vossieck list.

We start with giving a procedure for passing from the quiver with relations
for a tame concealed algebra $\L= \End_H(T)^{\op}$ to the quiver
of the cluster-tilted algebra $\G= \End_{\C_H}(T)^{\op}$, in our earlier notation.
The Happel-Vossieck list
gives the quivers together with a defining finite set of 
relations for the tame concealed algebras.
All of these relations have the following in common: they are of the form 
$\sum_{t=1}^r \sigma_i$ where $1 \leq r \leq 3$ and $\sigma_1, \dots, \sigma_r$ are
{\bf all} the paths from a vertex $i$ to a vertex $j$ and such that any
vertex different from $i,j$ is a vertex on at most one of the paths $\sigma_t$.
It is clear that each relation of this form is minimal. Recall that a relation $\rho$
is {\em minimal} if
whenever $\rho = \sum_{t=1}^n \alpha_t \rho_t \beta_t$, where
$\rho_1, \dots, \rho_n$ are relations, then for some $t$, both $\alpha_t$
and $\beta_t$ are scalars. We then have the following.

\begin{prop}\label{prop4.1}
Let $\L = \End_{H}(T)^{\op}$ be a tame concealed algebra, given as a quiver $Q$ and a
set of defining relations $\{ \rho_t \}$ from the Happel-Vossieck list. 
Then the quiver of the cluster-tilted algebra
$\G = \End_{\C}(T)^{\op}$ is obtained from $Q$ by adding an arrow
from the vertex $j$ to the vertex $i$   
if and only if one of the defining relations $\rho_t$ involves paths from $i$ to $j$. 
\end{prop}

\begin{proof}
Let $\rho = \rho_t = \sigma_1 + \cdots + \sigma_r$ be one of the defining relations in $\{ \rho_t \}$,
where $\{\sigma_1, \dots, \sigma_r \}$ are all paths from the vertex $i$ to the vertex $j$.

We first want to show that there is an arrow from $j$ to $i$ in the quiver of
$\G$.   
By reduction to finite type there is no arrow from a vertex $u$ to a vertex $v$ on $\sigma_1$,
where $(u,v) \neq (i,j)$, in the direction of $\sigma_1$. If there was an arrow $\beta$ in
the opposite direction, from $v$ to $u$ with $(u,v) \neq (i,j)$, we could choose $\beta$ such that
the cycle $u \to \cdots \to v \overset{\beta}{\to} u$ is full. Since $\sigma_1' \colon u \to \cdots \to v$
is not zero, but is zero in the factor algebra whose vertices are those of the path $\sigma_1'$, there must be 
another path $\sigma_1''$ from $u$ to $v$ in the quiver of $\G$, not going through a cycle. But then the
path from $i$ to $j$ having $\sigma_1''$ as a subpath would already have been part of $\rho$,
which gives a contradiction since $\sigma_1'$ and $\sigma_1''$ have more vertices in common than $i$ and $j$. 
Taking the factor algebra, keeping just the vertices of $\sigma_1$, it follows that $\sigma_1$ is a
minimal zero-relation for this factor, and hence there is an arrow from $j$ to $i$ in the 
quiver of $\G$. There can not be more than one arrow, since clearly if $\L$ is not
hereditary, it has at least 3 vertices, and since $\G$ is of minimal infinite type, there can be no double arrows 
in the quiver. 
There is no arrow from $i$ to $j$ since there is an arrow from $j$ to $i$, see \cite{bmr2}.

We want to show that there are no additional arrows. Let $\gamma \colon j \to i$ be
one of the arrows created from a relation, and let $\sigma_1, \dots , \sigma_r$, with $r \leq 3$, be the
corresponding shortest paths from $i$ to $j$ coming from the quiver of $\L$. Assume there is another 
shortest path $\psi$ from $i$ to $j$, which contains an arrow not in the quiver of $\L$.
By taking the factor keeping only the vertices of $\psi$, we see that $\psi$ is a minimal zero-relation in this 
factor, and hence $\psi$ is involved in a relation $\rho'$ from $i$ to $j$. If $\rho'$
was not minimal for $\G$, we would have $\psi = \alpha_1 \psi_1 \beta_1$
where $\psi_1$ is part of a relation from $u$ to $v$ with $(u,v) \neq (i,j)$. Taking
the factor corresponding to $\psi$ again, then $\psi_1$ would be zero, which is a contradiction.
Hence $\rho'$ is a minimal relation for $\G$. So we would, using shortening of paths as described in Section \ref{sec1.3},
have a cluster-tilted algebra with quiver
$$
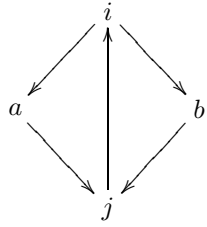
\xy
\xymatrix{
& i \ar[dl] \ar[dr] & \\
a \ar[dr] & & b \ar[dl] \\
& j \ar[uu] & 
}
\endxy
$$
with $i \to a \to j$ and $i \to b \to j$ being 0.
This is easily seen to be impossible.

Assume now that there is an arrow $\beta \colon i \to j$ in the quiver of $\G$ not
coming from a relation for $\L$. Since the quiver of $\L$ is connected,
there is a walk between $i$ and $j$ using only arrows from $\L$. Choose this walk as short as possible,
and adjust if necessary the choice of $\beta$ to make this walk shortest possible.
Consider the induced cycle,
and assume that it is not full. 
Then choose an arrow $\gamma$ different from $\beta$, such that the cycle created 
by $\beta$ and $\gamma$ is full. By the minimality of our choices, it follows
that $\gamma$ must be an arrow created from a relation for $\L$. Since by reduction to 
finite type this cycle must be oriented, this is impossible by the previous argument.
Hence the first cycle is full, and since it can not be the whole quiver $Q$,
it must be oriented.
Then for the corresponding 
factor there is a minimal zero-relation
$\sigma_1$ from $j$ to $i$. There is then a corresponding relation
$\rho = \sigma_1 + \cdots + \sigma_r$ for $\G$, from $i$ to $j$ (with as few paths as possible).
We claim that it is minimal.
If not, we have $\rho = \alpha_1 \rho_1 \beta_1 + \cdots + \alpha_s \rho_s \beta_s$, where
$\rho_1, \cdots , \rho_s$ are minimal relations and for each $t$, we have that $\alpha_t$ or $\beta_t$
is non-scalar. Then $\sigma_1$ must occur on the right hand side, in one of the relations, as some  
$\alpha_t \psi_t \beta_t$. Taking the factor corresponding to the vertices of $\sigma_1$, 
we get that $\psi_t$ is a zero-relation, contradicting the minimality of $\sigma_1$ 
in the factor. Hence $\beta$ is one of the arrows coming from a relation.
\end{proof}

The quivers of the cluster-tilted algebras $\End_{\C_H}(T)^{\op}$ are
the quivers of the minimal infinite cluster-tilted algebras by Theorem \ref{thm3.2},
which coincide with the minimal infinite cluster quivers by Theorem \ref{thm3.1}.
A classification of these is what A. Seven gave in \cite{se}.
Hence we get the following consequence.

\begin{thm}
By starting with the Happel-Vossieck list and replacing the given defining
relations with arrows in the opposite direction, we get the minimal infinite 
cluster quivers and hence the list of A. Seven.
\end{thm}

\section{Properties of minimal infinite tilted and cluster-tilted algebras}\label{sec-five}

Our next goal it to show which arrows we have to remove from a minimal
infinite cluster quiver in order to obtain the quiver with relations
for a minimal infinite tilted algebra. 
It turns out that this set is uniquely defined for each minimal infinite cluster quiver.

In this
section we give some necessary conditions on the set of arrows which
should be removed from the quiver of the minimal infinite 
cluster-tilted algebra. We also
discuss quadratic forms associated with signed graphs, and their
relationship to Tits forms in our context. Then the goal will be
completed in the next section, by proving a result on quadratic forms of 
signed graphs.

Let $T$ be a preprojective tilting module over a tame hereditary algebra
$H$. We investigate the passage from the tilted algebra $\L
=\End_{H}(T)^{\op}$ to the corresponding cluster-tilted
algebra $\G =\End_{\C_{H}}(T)^{\op}$ more carefully.
In particular, we would like information saying which arrows appearing in
the quiver $Q$ of $\Gamma$ are new ones.

\begin{prop}
With the above terminology, let $S$ be the set of new arrows
obtained when passing from $\Lambda$ to $\G$. Then $S$ contains
exactly one arrow from each full oriented cycle of $Q$, and no other
arrows.
\end{prop}

\begin{proof}
For each new arrow $\alpha$ we have exactly one, two or three 
full oriented cycles on which $\alpha$ lies.
All arrows except $\alpha$
on these cycles occur in the quiver of $\Lambda$, and there are no
other full cycles.
\end{proof}

It will be useful to look at quadratic forms associated 
with minimal infinite cluster quivers. We consider here quadratic forms
in $n$ variables $x_1, \dots ,x_n$ of the form
$q(x_1, \cdots ,x_n)= \sum_{i=1}^n x_i^2 +\sum_{i,j}a_{ij} x_ix_j
(i<j)$.
If $q$ is positive semidefinite, then the elements $z$ of $\ZZ^n$ with $q(z) = 0$
are called {\em radical vectors}.

For a finite-dimensional algebra $\L$ with quiver $Q$ without oriented cycles and with vertices
$1,\dots,n$ the corresponding Tits form is given by $a_{ij}=-t$ if $t$
is the number of arrows between $i$ and $j$ and $a_{ij}=s$ if $s$ is
the dimension of the space of minimal relations from $i$ to $j$ or
from $j$ to $i$. The Tits form of a tilted algebra of an extended
Dynkin quiver is known to be isomorphic to the one given by $Q$, and
is in particular positive semi-definite (see \cite{r}).

A crucial property of a tame concealed algebra is that there is a
positive sincere radical vector for the Tits form,
and amongst the tilted algebras of extended Dynkin type, the tame
concealed ones are exactly the ones with this property (see \cite{r}). 
The corresponding coordinates are given in the Happel-Vossieck list. 

Associated with any quadratic form there is a signed graph, i.e.
a graph with two kinds of edges which we call {\em solid} or {\em dotted}: If
$a_{ij}<0$ we have $-a_{ij}$ solid edges between $i$ and $j$, and if
$a_{ij}>0$ we have $a_{ij}$ dotted edges between $i$ and
$j$. Conversely, there is a quadratic form associated with a
signed graph in an obvious way (see \cite{r}). 

When passing from the quiver of a
tame concealed algebra $\Lambda$ to the quiver $Q$ of the corresponding
cluster-tilted algebra $\Gamma$,
with $S$ denoting the set of additional arrows, we associate with $Q$ and $S$
the signed graph obtained from $Q$ by making exactly the edges in $S$ dotted.
We denote by $q_S$ the associated quadratic form, which then clearly coincides with the Tits form
$t_{\L}$. Hence we have the following useful information.

\begin{prop}\label{pos-sin}
Let $\L = \End_{H}(T)^{\op}$ be a tame concealed algebra, and $Q$ the quiver of the corresponding cluster-tilted
algebra $\End_{\C}(T)^{\op}$, with $S$ the additional set of arrows. Then the quadratic form $q_S$ 
is isomorphic to the quadratic form of an extended Dynkin quiver mutation equivalent to $Q$ and
has a sincere
positive radical vector.
\end{prop}

Motivated by Proposition \ref{pos-sin} we call a set of arrows of a minimal 
infinite cluster quiver {\em admissible} if $S$ contains exactly one arrow from
each full oriented cycle, and no other arrows.
The strategy is to consider the quadratic form $q_S$
associated with $Q$ and $S$, and show that only one choice of admissible set will 
give a positive sincere radical vector. Here one is using that one can
show that $q_S$ is positive semidefinite, as we do in Section 6.

We denote by $\L_S$ the algebra whose quiver is obtained by removing the arrows in $S$
from $Q$, and with relations given as follows: For each arrow $\beta \colon j \to i$ in $S$, 
consider the sum of all shortest paths from $i$ to $j$. This coincides with the 
description of the tame concealed algebras in terms of quivers with relations for the 
``correct'' choice of $S$.
The only oriented cycles for $Q$, not containing properly an oriented cycle, are those 
created by the arrows $\alpha \colon j \to i$ corresponding to relations
on paths from $i$ to $j$. Since all paths from $i$ to $j$ not going through any cycle
are involved in the original relation,
it follows that all oriented cycles of the above type are full. Hence the quiver of $\L_S$
has no oriented cycles. Consider the relation $\rho$ associated with $\beta$.
We want to show that it is minimal for $\L_S$. Assume that we have 
$\rho = \sum_{t=1}^r \alpha_t \rho_t \beta_t$, where for each $t$ we have that 
$\rho_t$ is a relation associated with an arrow in $S$
and either $\alpha_t$ or $\beta_t$ is a non-trivial path. Let $\sigma_1$
be a path occurring for $\rho$. Then for some $t$, we have $\sigma_1 = \alpha_i \psi_i \beta_i$
where $\psi_i$ is a path in $\rho_i$. Since we have an arrow
from $v$ to $u$, where $\psi_i$ starts in $u$ and ends in $v$,
where $(u,v) \neq (i,j)$, we get a contradiction to $\sigma_1$ being a shortest path. 
Hence we have the following.

\begin{prop}
Let $Q$ be a minimal infinite cluster quiver, $S$ an admissible set of arrows and $\L_S$
the associated algebra. Then the quiver of $\L_S$ has no oriented cycles and the quadratic forms
$q_S$ and $t_S$ coincide.
\end{prop}

So dealing with the radical vectors of $q_S$ for an admissible set $S$ corresponds to dealing
with the
radical vectors of the Tits form $t_S = t_{\L_S}$ of the algebra $\L_S$. We shall in the next section
work with a more general choice of arrows $S$. 

\section{From cluster-tilted to tilted algebras}

Let $Q$ be a minimal infinite cluster quiver, with associated cluster-tilted
algebra $\G$. Denote by $Q_0$ the set of vertices of $Q$.
The aim of this section is to show the following:
For the quiver $Q$ there is a unique choice of an admissible set of
arrows $S$ such that the Tits form $t_S$ of the induced algebra
$\L_S$ has a positive sincere radical vector, or equivalently, such that
the quadratic form $q_S$ of the signed graph of $\G$ associated
with $S$ has a positive sincere radical vector.
An algorithm for finding this unique set $S$ is also given.
The result is seen as a consequence of a more general result  
about quadratic forms associated with signed graphs coming from 
minimal infinite cluster quivers.

The following sign change operation
at a vertex allows us to obtain a large number of isomorphic
quadratic forms on an undirected graph. Recall that two 
quadratic forms $q$ and $q'$ on
$\ZZ^{Q_0}$ are called isomorphic if there is an isomorphism $A$ on
$\ZZ^{Q_0}$ such that $q(A(v))=q'(v)$ for any $v$ in $\ZZ^{Q_0}$.

Let $\Sigma$ be a signed graph and let $i$ be a vertex in $\Sigma$.
We denote by $r_i(\Sigma)$ the graph obtained from $\Sigma$ by
changing the signs of the edges connected to the vertex $i$. If
$\Sigma = \Sigma(q)$ is the graph of a quadratic form $q$, then we denote by
$r_i(q)$ the quadratic form whose graph is $r_i(\Sigma)$.

Let us note the following obvious property of the sign change
operation.

\begin{prop}\label{pr:ref-parity}
Suppose that $C$ is a full cycle in the (undirected) graph
$\Sigma$. Then the parity of the number of dotted edges in $C$ is
the same both in $\Sigma$ and $r_i(\Sigma)$ for any vertex $i$ in
$\Sigma$ (thus $r_i$ preserves the parity of the number of dotted
edges in any cycle).
\end{prop}

We will mostly be interested in quadratic forms whose (signed)
graphs have the same underlying undirected unsigned graph, such as
those that can be obtained from each other by sign changes. It will
be convenient to use the following notation.
Let $Q$ be a quiver and $S$ a set of arrows of $Q$. We denote by
$q_S$ the quadratic form whose graph $\Sigma(q_S)$ is the underlying
(undirected) graph of $Q$ with the following sign assignment: any
edge whose corresponding arrow belongs to $S$ is dotted, and the
rest of the edges are solid.

Our next result characterizes a class of signed graphs that can be
obtained from each other by a special sequence of signs changes. It
is the main technical lemma that we use to prove the main theorem in this section.
Note that if a minimal infinite cluster quiver contains a full non-oriented
cycle, it must be of the form $\widetilde{A_n}$, and it comes from the same tame concealed algebra,
so that in this case there is nothing to prove.

\begin{lem}\label{lem:S-HV-2}
Suppose that $Q$ is a simply-laced cluster quiver which does not contain any non-oriented
full cycle. Let $S$ and $S'$ be two different sets of
arrows of $Q$ with the following property: for any full
cycle $C$, the parity of the number of arrows of $C$ contained in
$S$ is the same as the parity of the number of arrows contained in $S'$. Then the graph $\Sigma(q_{S'})$ of the
quadratic form $q_{S'}$ can be obtained from the graph $\Sigma(q_S)$
of $q_S$ by a sequence of sign changes such that a vertex is used at
most once and not all vertices are used. 
\end{lem}

\begin{proof}

We prove this by induction on the number, say $n$, of vertices of
$Q$:

For $n=2$, the underlying graph of the quiver $Q$ has a single edge,
which could be either solid or dotted, so we may assume that $q_S$
corresponds to the solid one and $q_{S'}$ to the dotted one. It is
obvious that sign change at any vertex transforms $\Sigma(q_S)$ to
$\Sigma(q_{S'})$.

Let us now assume that the lemma holds for quivers with $n-1$
vertices or less. Suppose that $Q$ has $n$
vertices.
Let us
consider a connected subquiver $\Sigma$ obtained by removing a
vertex, say $j$, from $Q$ (the existence of such a vertex leaving a
connected subquiver is easily seen). Since $\Sigma$ has
less than $n$ vertices, by the induction argument we have the
following: the restriction of $q_S$ to $\Sigma$ can be transformed
to the restriction of $q_{S'}$ to $\Sigma$ as described in the
lemma. Thus there is a sequence of mutually different vertices
$i_1,...,i_m$ in $\Sigma$ with $m<n-1$ such that
$$q_{S'}\mid_{\Sigma}= r_{i_m}...r_{i_1}(q_{S}\mid_{\Sigma}).$$ We claim
that either
\begin{equation}
\label{eq:lem1} q_{S'}= r_{i_m}...r_{i_1}(q_S),
\end{equation}
or
\begin{equation}
\label{eq:lem2} q_{S'}= r_{j}r_{i_m}...r_{i_1}(q_S).
\end{equation}
We need to show that if \eqref{eq:lem1} does not hold, then 
\eqref{eq:lem2} does.
So let us assume that
$q_{S'}$ is not equal to
$r_{i_m}...r_{i_1}(q_S)$. Denote this last expression by
$q_m$, i.e. $q_m=r_{i_m}...r_{i_1}(q_S)$. Note that, by our
induction assumption, the forms $q_S'$ and $q_m$ agree on the
subquiver $\Sigma$, i.e. (the signs of) the graphs of $q_{S'}$ and
$q_m$ on $\Sigma$ are the same. Since we assume $q_{S'}$ is not
equal to $q_m$, there is an edge $e$ (in the underlying graph of
$Q$) such that $e$ has opposite signs in $\Sigma(q_{S'})$ and
$\Sigma(q_m)$. Since $q_{S'}$ and $q_m$ agree on the subquiver
$\Sigma$ (obtained by removing the vertex $j$), the vertex $j$
must be one of the end points of $e$. We will show that, like $e$,
all of the edges containing $j$ have opposite signs in $q_{S'}$ and
$q_m$ (thus $q_{S'}= r_{j}(q_m)$, i.e. \eqref{eq:lem2} holds). For
this let us assume that there is an edge $f$, connected to $j$
such that $f$ has the same signs in $q_{S'}$ and $q_m$. Let us also
denote the remaining vertices of $e$ and $f$ by $r$ and $s$
respectively. We may also
assume that $e,f$ are such that the
vertices $r,s$ are such that the length of a shortest
possible walk is minimal.
Let us 
denote by $P$ a shortest walk connecting $r$ and $s$ in $\Sigma$.
We note that $j$ is not connected to any vertex on $P$ other than
$r$ and $s$ because of our assumption.
Thus the graph,
say $P_j$, induced by the path $P$ and the vertex $j$ is a
full cycle (which is oriented in $Q$ because we assumed that all
full cycles in $Q$ are oriented. Note also that $e$ and $f$ are the only
edges in $P_j$ which are not contained in $\Sigma$). We also note
that the parity of the dotted edges in $P_j$ is different in $q_{S'}$
and $q_m$; the difference is because $e$ has different signs in
$q_{S'}$ and $q_m$ and any other edge in $P_j$ has the same sign.
This gives a contradiction because sign change at a vertex preserves
the parity of dotted edges in cycles
(Proposition~\ref{pr:ref-parity}), thus $P_j$ has the same parity of
dotted edges in $q_m(=r_{i_m}...r_{i_1}(q_S))$ as in $q_S$, thus as
in $q_{S'}$, by our assumption. This completes the proof of the lemma. 
\end{proof}

Let us now give an algebraic description of the sign change
operation:

\begin{prop}\label{pr:reflection} 
Suppose that $q$ is a quadratic form.
Then $r_k(\Sigma(q))$ is the graph of the form $q$
with respect to the basis (variables) obtained from the basis
(variables) for $\Sigma(q)$ by changing $x_k$ to $-x_k$ and keeping
the other elements the same.
\end{prop}

\begin{proof} 
Let us assume without loss of generality that $k=1$. Let
$\{y_i\}$ be the dual basis such that $y_1=-x_1$ and $y_i=x_i$ for
$i \ne 1$. We note that we have
$$q=a_{1,1}x_1^2 +\sum a_{1,i}(x_1) x_i + \sum a_{i,j}x_i x_j, 2\leq
i \leq j,$$ or equivalently
$$q=a_{1,1}(-x_1)^2 +\sum -a_{1,i}(-x_1)
x_i + \sum a_{i,j}x_i x_j, 2\leq i \leq j$$ thus
$$q=a_{1,1}y_1^2 +\sum -a_{1,i}y_1
x_i + \sum a_{i,j}y_i y_j, 2\leq i \leq j.$$ Then the statement
follows from our definitions.
\end{proof}

\begin{cor}\label{cor:reflection}
If $q$ is positive semi-definite of corank 1, then the coordinates
of the radical vector of $q$ with respect to the basis corresponding
to the graph $r_i(\Sigma(q))$ is the same as the one for $\Sigma(q)$
with the exception that the $i-th$ coordinate is the negative of the
one for $\Sigma(q)$.
\end{cor}

The following basic fact gives interesting examples of signed graphs where
each full cycle has an odd number of dotted edges.

\begin{prop}\label{pr:odd-positive}
If $q$ is a positive-definite form, then any full cycle 
in $\Sigma(q)$ has an odd number of dotted edges.
\end{prop}

\begin{proof} 
Suppose that $\Sigma(q)$ contains a full cycle $C$ which has
exactly an even number of dotted edges. To prove the proposition, it
is enough to show that the restriction of $q$ to $C$ is positive
semi-definite. By Lemma~\ref{lem:S-HV-2}, the form $q$ is isomorphic
to the form $q'$ such that the underlying (undirected and unsigned)
graph of $\Sigma(q' \mid_C)$ is the same as $C$ with no dotted edges.
Let $u$ be the vector such that for any vertex of $C$ the
corresponding coordinate of $u$ is $1$ and the rest are $0$. Then
$q'(u)=0$, so $q'$ is not positive definite, thus $q$ is not
positive definite either.
\end{proof}

The classification of positive-definite quadratic
forms is the well-known Cartan-Killing classification: any
positive-definite quadratic form is isomorphic to one given
by a Dynkin graph. Similarly positive-semidefinite 
quadratic forms of corank 1 are classified by the extended Dynkin
graphs \cite{r}.

If $Q$ is a finite type cluster quiver, then the form $q_S$ is
positive-definite. More precisely, the following follows from \cite[Thm 1.1]{bgz}.

\begin{prop}\label{pr:finite-positive}
Suppose that $Q$ is a finite type cluster quiver (note that $Q$ does not
have any non-oriented full cycles) and let $S$ be a set of arrows of $Q$
such that $S$ contains exactly an odd number of arrows from each
oriented cycle. Then the
form $q_S$ is positive definite. 
\end{prop}

We now show the main result of this section.

\begin{thm}\label{th:S-HV-2}
Let $Q$ be a minimal infinite cluster quiver. For any set $S$
of arrows such that $S$ contains exactly an odd number of arrows
from each oriented cycle and no arrows from any non-oriented cycle,
the corresponding quadratic form $q_S$ is isomorphic to the one
defined by an extended Dynkin quiver which is mutation equivalent to
$Q$. Furthermore, there is a unique set $S_+$ where the form
$q_{S_+}$ has a sincere positive radical vector.
\end{thm}

\begin{proof}
We can assume that $Q$ has at least three vertices, since for
two vertices the set $S$ would be empty. 
It follows from Proposition \ref{pos-sin} and Lemma \ref{lem:S-HV-2} that for any $S$
as in the theorem, $q_S$ is isomorphic to the quadratic form given by an extended 
Dynkin quiver mutation equivalent to $Q$. Let $S_+$ be a set such that 
$\L_{S_+}$ is tame concealed, hence $q_{S_+}$ has a sincere positive radical vector.
Now take $S'$ as in the theorem such that
$S'\ne S_+$. Then there is a sequence of vertices $i_1,...,i_m$ as
in Lemma~\ref{lem:S-HV-2} such that $q_{S'}=
r_{i_m}...r_{i_1}(q_{S_+})$, ($m\geq 1$). Since a sign change at a
vertex changes only the corresponding coordinate of the radical
vector to its negative (Corollary~\ref{cor:reflection}) and each
vertex is used once (and not all vertices are used), exactly the
coordinates corresponding to the vertices $i_1,...,i_m$ will be
negative in the radical vector of $q_{S'}$ (and the remaining ones
will be positive). This completes the proof.
\end{proof}

In view of the discussion in Section \ref{sec-five},
we have the following reformulation of Theorem \ref{th:S-HV-2}.

\begin{thm}
Let $Q$ be a minimal infinite cluster quiver.
Then, for any admissible set of arrows $S$, the Tits form $t_S$
of the algebra $\L_S$ is isomorphic to the quadratic form 
defined by the extended Dynkin quiver which is mutation equivalent to $Q$.
Furthermore, there is a unique admissible set $S_+$ such that $t_{S_+}$
has a positive sincere radical vector. Also, for any
admissible $S$, the algebra $\L_S$ is of minimal infinite representation type if
and only if $S = S_+$.
\end{thm}

Note that if we are given a minimal infinite cluster quiver 
$Q$, and we find a set of arrows $S$, with exactly one arrow from each oriented cycle
and no other arrows, such that the quadratic form of the associated signed graph has a positive sincere
radical vector, then we know by Theorem \ref{th:S-HV-2} that this is the correct choice 
for obtaining the quiver of the associated tame concealed algebra.

More systematically, as a consequence of the proof
of Theorem \ref{th:S-HV-2}, we have the following simple procedure for finding the 
desired set $S_+$. 

\begin{algorithm}\label{alg:} 
To find $S_+$ in a minimal infinite cluster quiver $Q$:

\begin{itemize}
\item[(1)]
Take any set of arrows $S$ which contains exactly an odd number of
arrows from any oriented full cycle,
\item[(2)]
compute the radical vector of the quadratic form $q_S$,
\item[(3)]
apply sign change operation at all vertices whose corresponding
coordinates in the radical vector of $q_S$ is negative,
\item[(4)]
the set of arrows corresponding to the dotted edges is exactly
$S_+$.
\end{itemize}

\end{algorithm}

We end this section by an example and a conjecture. The following seems to be
a general fact: If we remove an admissible $S$ such that $\L_S$
is not of minimal infinite type, then $\L_S$ is of finite type. 
To give an example consider the following minimal infinite cluster quiver.
$$
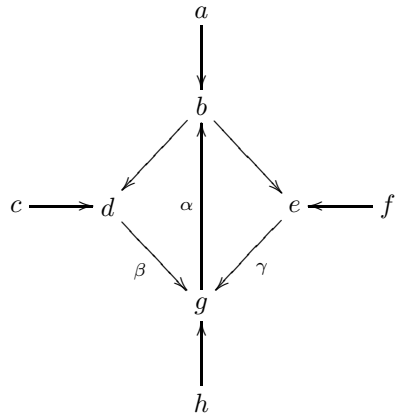
\xy
\xymatrix{
& & a \ar[d] & & \\
& & b \ar[dl] \ar[dr] &  &\\
c \ar[r] & d \ar_{\beta}[dr] & & e \ar^{\gamma}[dl] & f \ar[l] \\
& & g \ar^{\alpha}[uu] & & \\
& & h \ar[u] & & 
}
\endxy
$$
The choice of $S$ giving $\L_S$ of minimal infinite type is $S = \{\alpha\}$.
Another admissible set is $S' = \{ \beta, \gamma \}$. In this case $\L_{S'}$ 
is of global dimension 3, and hence not tilted, and of finite representation type.

We make the following conjecture.

\begin{conj}
In the set-up of Theorem \ref{th:S-HV-2}, for any admissible set $S$,
the dimension vectors of the indecomposable representations of $Q_S$
are exactly the positive roots of the corresponding Tits form $t_S$
(here $x$ is a root if $t_S(x)=1$). Furthermore if $S \neq S_+$, then the algebra
$\L_S$ is of finite representation type.
\end{conj}

\section{The cluster-tilted algebras of minimal infinite cluster quivers} 

In \cite{bmr3} it was shown that a cluster-tilted algebra of finite type
(over an algebraically closed field $k$) is uniquely determined by its quiver,
and has relations of a nice form. In this final section we show that
similar results hold for the cluster-tilted algebras investigated in this paper.
Recall that all minimal infinite cluster-tilted algebras are 
endomorphism algebras of preprojective tilting modules over tame hereditary algebras.

\begin{thm}\label{thm5.1}
Let $T$ be a preprojective tilting module over a tame
hereditary algebra $H$, and $\Gamma
=\End_{\mathcal{C}_{H}}(T)^{\op}$ the associated (basic) minimal
infinite cluster-tilted algebra. For each arrow $\alpha
\colon j \to i$ lying on an oriented cycle in the quiver of $\Gamma$, 
we consider all the shortest paths
$\sigma_1,\cdots ,\sigma_r$ from $i$ to $j$.
Then there is
a minimal relation $\sigma_1 + \cdots + \sigma_r$, and $r \leq 3$.
Furthermore, the set of relations obtained this way is a generating set of minimal relations.
\end{thm}

\begin{proof}

We will use that a path passing through an oriented cycle is a zero-path. To see this,
we use that such a path corresponds to an endomorphism $X \to X$ of a preprojective indecomposable module $X$.
In the derived category $\Hom_{D^b(H)}(X,FX)$ is clearly zero, so the claim holds.

We first show that for every arrow $\alpha \colon j \to i$ lying on a cycle, there is 
a relation of the prescribed form. 

In case there are two or more shortest paths from $i$ to $j$, we claim that the subquiver 
generated by these paths has the following shape ($\star$).
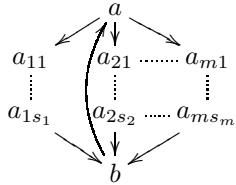
$$
\xy
\xymatrix@C0.3cm@R0.3cm{
  & a \ar[dl] \ar[d] \ar[dr] & \\
  a_{11} \ar@{..}[d] & a_{21} \ar@{..}[d] \ar@{..}[r] & a_{m1} \ar@{..}[d] \\
  a_{1s_1} \ar[dr] & a_{2s_2} \ar[d]  \ar@{..}[r] & a_{ms_m} \ar[dl]\\
  & b \ar@/^11pt/[uuu] &
}
\endxy
$$
To prove the claim, we apply Lemma 2.14 in \cite{bmr3} and note that the proof goes
through in our setting. 

We note that there are minimal infinite cluster quivers of this form with three ``arms''
from $a$ to $b$ (i.e. $m=3$), but that it is evident that $m \geq 4$ would give an algebra
which is not of minimal infinite type, since there would be too many arrows meeting in $a$,
i.e. after removing the vertex $b$, there would still four arrows all starting in $a$.

We first discuss quivers which are not of the form ($\star$).
In this case it is clear by the above that
for any arrow
$\alpha \colon j \to i$ lying on a cycle in $Q$, there is
a vertex $e$ such that $e$ does not lie
on any of the shortest paths from $i$ to $j$. 

Consider the factor algebra $\Gamma /{\Gamma e \Gamma}$.
This is a cluster-tilted algebra (by \cite{bmr2}) of finite representation type.
Thus, using the main result from \cite{bmr3} there is a  
minimal relation $\rho$ in the factor algebra involving the shortest paths from $i$ to $j$.
Also from \cite{bmr3} it follows that there are at most two such shortest paths. Note
that this can also be observed from the Happel-Vossieck list, or the A. Seven list. 
Also it is shown in \cite{bmr3} that if there are two such shortest paths $\rho_1$ and $\rho_2$,
then the corresponding relation is $\rho_1 - \rho_2$. It is clear from
inspection of the Happel-Vossieck list, that we get an isomorphic algebra by changing the
relation to $\rho_1 + \rho_2$.

We now claim that $\rho$ is also a minimal
relation for $\Gamma$. 
To prove the claim, 
we first observe that $\rho$ is a relation for $\Gamma$. For since
$\rho$ is a  relation for the factor algebra $\Gamma /{\Gamma e
\Gamma}$, there is a sum of paths $\rho'$ for $\Gamma$, all going
through the vertex $e$, such that $\rho +\rho'$ is a relation for
$\Gamma$. We claim that any such path 
from $i$ to $j$ which is not a shortest path must 
pass through an oriented cycle, thus it is itself a $0$-relation.
Consequently $\rho$ is a relation for
$\Gamma$. 

To prove this claim,
assume there is a path $\psi \colon i \to a_1 \to \cdots \to a_s = j$ which is not shortest and does
not pass through an oriented cycle. We can assume that the path has minimal length among all paths
from $i$ to $j$ with this property.
Consider the subquiver

$$
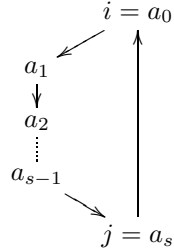
\xy
\xymatrix@C0.3cm@R0.3cm{
  & i = a_0  \ar[dl] \\
  a_1 \ar[d] & \\
  a_2 \ar@{..}[d] & \\
  a_{s-1} \ar[dr] & \\
  & j = a_s\ar[uuuu] &
}
\endxy
$$

If there are further arrows going downwards,
say $a_x \to a_y$ with $y> x+1$, this would contradict the minimality of the length 
of $\psi$, so no such arrows exist.
If there are 
additional arrows going upwards, choose such an arrow starting
in a vertex $a_y$ with $y$ maximal. It is clear that the quiver has a factor quiver
which is a non-oriented cycle, and can hence not be the quiver of a cluster-tilted algebra
of finite representation type. This is a contradiction, so all paths which are
not shortest pass through an oriented cycle.

To show that $\rho$ is minimal for $\Gamma$, assume $\rho
=\alpha_1 \rho_1 \beta_1 + \dots + \alpha_n \rho_n \beta_n$, where
$\rho_1, \cdots, \rho_n$ are minimal relations for $\Gamma$. If a
path going  through $e$ occurs for some $\alpha_i \rho_i \beta_i$,
then this path goes  through a cycle, and is
hence itself a $0$-relation, which is not minimal. Since no such
path occurs for $\rho$, these paths can be removed on the right hand
side. Hence for some $i$ we have that $\alpha_i$ and $\beta_i$ are
constant, using that $\rho$ is a minimal  relation for $\Gamma
/{\Gamma e \Gamma}$. Then we conclude that $\rho$ is minimal also for
$\Gamma$.

We now consider quivers of minimal infinite cluster-tilted algebras
which have a quiver of the form ($\star$).
Then for all arrows except the arrow $b \to a$, we can use the same technique as in the previous case,
to conclude that there is a minimal relation as prescribed.
For the arrow $b \to a$ it is clear that the corresponding tilted algebra is the
algebra with quiver obtained by removing the arrow $b \to a$, and there is a relation
as prescribed. 

Next we need to prove that all the minimal relations give rise to arrows. 
More precisely, let $\rho$ be a minimal relation involving paths from $i$ 
to $j$ in the cluster-tilted algebra. We need to show that there exists an arrow
from $j$ to $i$. For this the following lemma
due to Assem, Br{\"u}stle, Schiffler and Reiten, Todorov is useful, and gives a
simplification of our original proof. 

\begin{lem}
Let $A$ be any cluster-tilted algebra, and let $S_1, S_2$ be simple $A$-modules.
Then $\dim_k \Ext_A^1(S_1, S_2) \geq \dim_k \Ext_A^2(S_2,S_1)$. 
\end{lem} 

Recall that the number of arrows from the vertex corresponding to $S_1$ to the vertex corresponding to $S_2$
is given by $\dim_k \Ext_A^1(S_1, S_2)$ and that $ \dim_k \Ext_A^2(S_2,S_1)$
is the dimension of the space of minimal relations involving paths from
the vertex corresponding to $S_2$ to the vertex corresponding to $S_1$, by \cite{bo}.
Thus the proof of the theorem is finished by the above lemma.
\end{proof}

Remark: Actually, the last part of the proof of the theorem, namely 
that minimal relations give rise to arrows, can also be seen
directly from studying the quivers and relations on the Happel-Vossieck list. 
But for a few of the quivers on this list this is a rather cumbersome procedure. We therefore
included the more general lemma above. 

Using the description of the relations given in Theorem \ref{thm5.1}, it is easy to see that
the algebras $\L_S$, for an admissible set $S$, as defined in Section \ref{sec-five},
are obtained by restriction of the relations from the cluster-tilted algebra, in the case of 
minimal infinite type.

\section*{Acknowledgments}

This work started when the third author A. Seven
visited NTNU in spring 2005 as a LieGrits postdoctoral fellow.
He thanks the coauthors I. Reiten and A. Buan for their
kind hospitality. He also thanks the members of the
algebra group at NTNU for many helpful discussions.

\end{document}